\documentclass[11pt]{article}\synctex=1 
\usepackage{amsfonts}
\usepackage{amssymb}
\usepackage{amsmath}
\usepackage{amsthm}
\usepackage{geometry}
\usepackage{commath}
\usepackage{booktabs}
\usepackage{lscape}
\usepackage{xr}
\usepackage[authoryear]{natbib}

\newcommand{\myref}[2]{\hyperref[#1]{#2}}
\numberwithin{equation}{section}

\usepackage{latexsym}
\usepackage{graphicx} 
\usepackage{enumerate}
\usepackage{setspace}
\usepackage[toc,title,titletoc,header]{appendix}
\usepackage[bottom]{footmisc}   
\usepackage[pdfborder={0 0 0}]{hyperref}

\usepackage{lineno}
\setpagewiselinenumbers

\newtheorem{theorem}{Theorem}[section]

\theoremstyle{definition}

\theoremstyle{remark}

\newcounter{assumptionM}
\newcounter{assumptionA}
\def\theassumptionM{M.\arabic{assumptionM}}
\def\theassumptionA{A.\arabic{assumptionA}}
\hypersetup{
  colorlinks=true,linkcolor=blue, citecolor=blue, filecolor=blue, 
  linkbordercolor={0 1 1}, citebordercolor={1 0 0},
  pdftitle={Inference in partially identified models with many moment inequalities using Lasso},
  pdfauthor={Bugni, Caner, Kock, and Lahiri},
  pdfcreator={\LaTeX\ with package \flqq hyperref\frqq},
  pdfsubject={},
}

\begin{document}
\sloppy

\title{
Delta Theorem in the Age of High Dimensions\thanks{ Department of Economics, 452 Arps Hall, TDA  Columbus, OH, 43210. email:caner.12@osu.edu}
}

\author{
Mehmet Caner\\
Department of Economics\\
Ohio State University
}

\maketitle

\begin{abstract}
We provide a new version of delta theorem, that takes into account of high dimensional parameter estimation. We show that depending on the structure of the function, the limits of functions of estimators have faster or slower rate of convergence than the limits of estimators. 	We illustrate this via two examples. First, we use it  for testing in high dimensions, and second in estimating large portfolio risk. Our theorem works  in  the case of larger number of parameters, $p$, than the sample size, $n$: $p>n$.
	


\end{abstract}



\newpage
\section{Introduction}
Delta Method is one of the most widely used theorems in econometrics and statistics. It is a very simple and useful idea. It can provide limits for complicated functions of estimators as long as function is differentiable. Basically, the idea is the limit of the function of estimators can be obtained from the limit of the estimators, and with exactly the same rate of convergence. In the case of finite dimensional parameter estimation, due to  derivative at the parameter value being finite,  rates of convergence of both estimators and function of estimators are the same. 

In the case of high dimensional parameter estimation, we show that this is not the case, and the rates of convergence may change. We show that the structure of the function is the key, and depending on that functions of estimators  may converge faster or slower than estimators. In this paper, we provide a new version of delta method which accounts for high dimensions, and generalizes the previous finite dimensional case. After that we illustrate our point in two examples:  first by  examining a linear function of estimators that is heavily used in econometrics, and second by analyzing the risk of a large portfolio of assets in finance. Section 2 provides new delta method. Section 3 has two examples.  Appendix shows the proof.

\section{High Dimensional Delta Theorem}

Let $\beta_0 = (\beta_{10}, \cdots, \beta_{p0})'$ be a $p \times 1$ parameter vector with an estimator $\hat{\beta} = (\hat{\beta}_1, \cdots, \hat{\beta}_p)'$.  Define a function $f(.)$, $f: K \subset R^p \to R^m$ defined at least on a neighborhood of $\beta_0$, where $p>m$, $m$ will be taken as constant  for convenience, but $p$ is allowed to increase when $n$ increases. Furthermore the function $f(.)$ is differentiable  at $\beta_0$, which means, with $h \neq 0$, 
\[ \lim_{h \to 0} \frac{\| f(\beta_0 + h ) -  f(\beta_0) - f_d (\beta_0) h \|_2}{\|h \|_2} = o(1),\]
where $\|.\|_2$ is the Euclidean norm for a generic vector, and $f_d (\beta_0)$ is the $m \times p$ matrix, which $ij$th cell consists of $\partial f_i / \partial \beta_j$ evaluated at $\beta_{0}$, for $i=1,\cdots, m$, $j=1,\cdots p$. 

Before the theorem, we need the following matrix norm inequality. Take a generic matrix which is of dimension $m \times p$. Denote the Frobenius norm for a matrix 
as $\|| A\||_2 = \sqrt{\sum_{i=1}^m \sum_{j=1}^p a_{ij}^2}$. Note that in some of the literature such as Horn and Johnson (2013), this definition  is not considered a matrix norm, due to lack of submultiplicativity. However, our results will not change regardless of matrix norm definitions, if we abide by Horn and Johnson (2013), our results can be summarized in algebraic form, rather than matrix norm format. Define 
\[ A = \left[ \begin{array}{c}
		a_1' \\
		\vdots \\
		a_m' \end{array} \right],\]
		where $a_i$ is $p \times 1$ vector, and its transpose is $a_i'$, $i=1,\cdots,m$.
		Then for a generic $p \times 1$ vector x,
		\begin{equation}
		\| A x \|_2 = \sqrt{\sum_{i=1}^m (a_i' x)^2} \le \left( \sqrt{\sum_{i=1}^m \|a_i \|_2^2}\right) \|x\|_2 = \|| A \||_2 \|x\|_2,\label{1}
		\end{equation}
where the inequality is obtained by Cauchy-Schwarz inequality. Note that if we apply Horn and Johnson (2013) norm definition, this matrix norm inequality still holds, but we cannot use the matrix norm.  In that case we have 
\begin{equation}
\| A x \|_2 \le \left( \sqrt{\sum_{i=1}^m \|a_i \|_2^2}\right) \|x\|_2 .\label{2}\end{equation}

Our new delta theorem is provided for high dimensional case. This generalizes  Theorem 3.1 of van der Vaart (2000). Key element in our Theorem below is 
$\||f_d (\beta_0)\||_2$. We should note that this norm of matrix derivative depends  on $n$, through $p$, which is the number of columns in $f_d (\beta_0)$.
Let $r_n, r_n^* \to \infty$, as $n \to \infty$, be the rate of convergence of estimators, and the functions of estimators, respectively in the theorem below.

\begin{theorem} Let  a function $f(\beta): K \subset R^p \to R^m$, and differentiable at $\beta_0$. Let $\hat{\beta}$ be the estimators for $\beta_0$, and $\hat{\beta} \neq \beta_0$, assume we have the following result:
\[ r_n \| \hat{\beta} - \beta_0 \|_2 = O_p (1).\]

a) Then if  $\|| f_d (\beta_0) \||_2 \neq o(1)$, with $\|| f_d (\beta_0) \||_2 > 0$, we get 
\begin{equation}
 r_n^* \| f(\hat{\beta}) - f (\beta_0) \|_2 = O_p (1),\label{3}
 \end{equation}
where   
\begin{equation}
 r_n^* =  O \left(\frac{r_n}{  \|| f_d (\beta_0) \||_2} \right).\label{4}
 \end{equation}
 
 b) Then if $\|| f_d (\beta_0) \||_2 = o(1)$, we get
 \begin{equation}
 r_n \| f (\hat{\beta}) - f (\beta_0) \|_2 = o_p (1).\label{4b}
 \end{equation} 
\label{thm1} 
\end{theorem}
Remarks. 1. Note that in part a), with $r_n^*$, we have a slower or the same rate of convergence as in $r_n$. In part b), clearly, the function of estimators converge to zero in probability faster than the rate of estimators themselves. This can be seen from noting that in part b), even though 
\[ \| \hat{\beta} - \beta_0 \|_2 = O_p (1/r_n),\]
for functions
\[\| f(\hat{\beta}) - f(\beta_0) \|_2 = o_p (1/r_n).\]

2. Note that Horn and Johnson (2013) defines Frobenius norm only for square matrices unlike our case. If we use their approach, then our main result in part a) will be 
\begin{equation}
 r_n^*= O \left( \frac{r_n}{  \sqrt{\sum_{i=1}^m \| f_{di} (\beta_{0}) \|_2^2} }\right),\label{r1}
 \end{equation}
where we use (\ref{2}) instead of (\ref{1}) in the proof of Theorem \ref{thm1}a. Also note that $f_{di} (\beta_{0})$ is the $p \times 1$ vector, which is $\partial f_i(.)/ \partial \beta$ evaluated at $\beta_0$.

3. Also see that this result (\ref{4}) can be obtained  in other matrix norms subject to the same caveat in Remark 1. A simple Holder's inequality provides
\begin{equation}
\| A x \|_1 \le \|| A \||_1 \|x\|_1 ,\label{r2}\end{equation}
where  we define the maximum column sum matrix norm: $\||A\||_1 = \max_{1 \le j \le p} \sum_{i=1}^m | a_{ij}|$, where $a_{ij}$ is the $ij$ th element of $A$ matrix.
Applying this in the proof of Theorem \ref{thm1}a, given $r_n \| \hat{\beta} - \beta_0 \|_1 = O_p (1)$ and replacing everything with Frobenius norm for matrices with maximum column sum matrix norm, and $l_2$ norm for vectors with $l_1$ norm, we have 
\begin{equation}
r_n^* =  O \left(\frac{r_n}{  \|| f_d (\beta_0) \||_1}\right).\label{r3}
 \end{equation}
 Also part b) can be written in $l_1$ norm as well.
 
 4. We can also extend these results to another norm.  A simple inequality provides
\begin{equation}
\| A x \|_{\infty}  \le \|| A \||_{\infty} \|x\|_{\infty} ,\label{r4}\end{equation}
where  we define the maximum row sum matrix norm: $\||A\||_{\infty} = \max_{1 \le i \le m} \sum_{j=1}^p | a_{ij}|$, where $a_{ij}$ is the $ij$ th element of $A$ matrix.
Applying this in the proof of Theorem \ref{thm1}a, given $r_n \| \hat{\beta} - \beta_0 \|_{\infty}  = O_p (1)$ and replacing everything with Frobenius norm for matrices with maximum row sum matrix norm, and $l_2$ norm for vectors with $l_{\infty}$ norm, we have 
\begin{equation}
r_n^* = O \left(\frac{r_n}{  \|| f_d (\beta_0) \||_{\infty}}\right).\label{r4}
 \end{equation}
 Also part b) can be written in $l_{\infty}$ norm as well.
 
 5. What if we have $m=p$? or $m>p$, and $m \to \infty$ as $n \to \infty$? Then all our results will go through as well, this is clear from our proof.

\section{Examples}

We now provide two examples that will highlight the contribution. First one is related to linear functions of estimators in part a), and the second one is related to risk of the large portfolios, and part b).

{\bf Example 1}. 

Let us denote $\beta_0$ as the true value of vector ($p \times 1$) of coefficients. The number of the true nonzero coefficients are denoted by $s_0$, and $s_0>0$. 
A simple linear model is:
\[ y_t = x_t' \beta_0 + u_t,\]
where $t=1,\cdots, n$, with $u_t$ iid mean zero, finite variance error, and $x_t$ is deterministic set of $p$ regressors for ease of analysis.

The lasso estimator in a simple linear model is defined as
\[ \hat{\beta} = argmin_{\beta \in R^p} \sum_{t=1}^n \frac{(y_t - x_t' \beta)^2}{n} + \frac{\lambda}{n} \sum_{j=1}^p | \beta_{j0}|,\]
where $\lambda$ is a positive tuning parameter, and it is established that $\lambda = O (\sqrt{\frac{logp}{n}})$.  
Corollary 6.14 or Lemma 6.10 of Buhlmann and van de Geer (2011) shows, for lasso estimators $\hat{\beta}$, with $p>n$
\begin{equation}
r_n \| \hat{\beta} - \beta_0 \|_2 = O_p (1),\label{ex11}
\end{equation}
where 
\begin{equation}
r_n = \sqrt{\frac{n}{logp}} \frac{1}{\sqrt{s_0}}.\label{ex11a}
\end{equation}
At this point, we will not go into detail such as what assumptions are needed to get (\ref{ex11}), except to tell that minimal adaptive restrictive eigenvalue is positive, and noise reduction is achieved. Details can be seen in Chapter 6 in Buhlmann and van  de Geer (2011). 

The issue is what if the researchers are interested in the asymptotics of $D (\hat{\beta} - \beta_0)$, where $D: m \times p$ matrix. $D$ matrix can be thought of putting restrictions on $\beta_0$. We want to see whether $D (\hat{\beta} - \beta_0)$ has a different rate of convergence than $\hat{\beta} - \beta_0$. From our Theorem \ref{thm1}a, it is clear that $f_d (\beta_0) = D$. Assume that $\||D\||_2 \neq o(1)$. So 
\[ r_n^* \| D (\hat{\beta} - \beta_0) \| = O_p (1),\]
where 
\begin{equation}
r_n^* = O \left(\frac{r_n}{\|| D \||_2} \right).\label{ex12}
\end{equation}
We know from matrix norm definition:
\[ \|| D \||_2 = \sqrt{\sum_{i=1}^m  \| d_i \|_2^2},\]
and $\|d_i\|_2$ is the Euclidean norm for vector $d_i$ which is $p \times 1$, and 
\[ D = \left[ \begin{array}{c}
                    d_1' \\
                    \vdots \\
                    d_m' \end{array} \right].\]
                    Basically in the case of inference, this matrix and vectors show how many of $\beta_0$ will be involved with restrictions. If we want to use $s_0$ elements in each row of $D$  to test $m$ restrictions, then $ \|| D \||_2 = O (\sqrt{s_0})$.  Note that this corresponds to using $s_0$ elements in $\beta_0$ for testing $m$ restrictions. 
                    So 
                    \begin{equation}
                    r_n^* = \left(\sqrt{\frac{n}{logp}} \right) \left( \frac{1}{s_0}\right),\label{ex13}
                    \end{equation}
                    which shows that even though we have fixed number of restrictions, $m$, using $s_0$ of coefficients in testing will slow down the rate of convergence, $r_n^*$  by $\sqrt{s_0}$, compared with lasso estimators, rate of $r_n$.  This can be seen by comparing (\ref{ex11a}) with (\ref{ex13}).

{\bf Example 2}. 

One of the cornerstones of the portfolio optimization is estimation of risk. If we denote the portfolio allocation vector by $w$ ($p \times 1$) vector, and the covariance matrix of asset returns by $\Sigma$, the risk is $\sqrt{(w' \Sigma w)}$. We want to analyze risk estimation error which is $(w' \hat{\Sigma} w)^{1/2} - 
(w' \Sigma w)^{1/2}$.  $\hat{\Sigma}$ is the sample covariance matrix of asset returns. We could have analyzed  risk estimation error with estimated weights, as in Fan etal (2015), 
$(\hat{w} \hat{\Sigma} \hat{w})^{1/2} - (\hat{w} \Sigma \hat{w})^{1/2}$, but this extends the analysis with more notation with the same results.

A crucial step  in assessing the accuracy of risk estimator is given in p.367 of Fan etal (2015), which is  the term
$w' (\hat{\Sigma} - \Sigma) w$. . Just to simplify the analysis, we will assume iid, sub-Gaussian asset returns. Also we will find the global minimum variance portfolio as in Example 3.1 of Fan etal (2015). So
\[ w = \frac{\Sigma^{-1} 1_p}{1_p' \Sigma^{-1} 1_p},\]
where $\Sigma$ is nonsingular, and $1_p$ is the $p$ vector of ones. Assume $0 < Eigmin (\Sigma^{-1}) \le  Eigmax (\Sigma^{-1}) < \infty$, where $Eigmin(.), Eigmax(.)$ represents the minimal, maximal eigenvalues respectively of the matrix inside the parentheses. In Remark 3 of Theorem 3 
in Caner etal (2016) 
\begin{equation}
\| w \|_1 = O ( \max_j  \sqrt{s_j}),\label{ex21}
\end{equation}
where $s_j$ is the number of nonzero cells in $j$ th row of $\Sigma^{-1}$ matrix, $j=1,\cdots, p$. Equation (\ref{ex21}) represents case of growing exposure, which means we allow for extreme positions in our portfolio, since we allow $s_j \to \infty$, as $n \to \infty$. See that 
\begin{equation}
\| \hat{\Sigma} - \Sigma \|_{\infty} = O_p (\sqrt{\frac{logp}{n}}),\label{ex22}
\end{equation}
by van de Geer etal. (2014). Then clearly by (\ref{ex21})(\ref{ex22})

\begin{equation}
| w' (\hat{\Sigma}- \Sigma) w | \le \| w\|_1^2 \|\hat{\Sigma} - \Sigma \|_{\infty} = O_p (\max_{1 \le j \le p} s_j \sqrt{logp/n}).\label{ex23}
\end{equation}
This means taking $\beta_0= w' \Sigma w, \hat{\beta} = w' \hat{\Sigma} w$, so $m=1$  in  Theorem \ref{thm1}b,
\begin{equation}
r_n | w' (\hat{\Sigma}  - \Sigma ) w| = O_p (1),\label{ex24}
\end{equation}
where 
\begin{equation}
r_n = \sqrt{\frac{n}{logp}} \frac{1}{\max_{1 \le j \le p} s_j}.\label{ex25}
\end{equation}

But the main issue is to get risk estimation error, not the quantity in (\ref{ex24}). To go in that direction see that 
\begin{equation}
O (\max_{1 \le j \le p} s_j)= \|w\|_1^2 Eigmin(\Sigma) \le | w' \Sigma w| \le \|w\|_1^2 Eigmax( \Sigma) = O (\max_{1 \le j \le p} s_j),\label{ex26}
\end{equation}
where we use (\ref{ex21}) and $Eigmax (\Sigma) < \infty$, $Eigmin(\Sigma)>0$.   

Note   that risk is $f(\beta_0) = (w' \Sigma w)^{1/2}$, and $f_d (\beta_0)=(w' \Sigma w)^{-1/2}= O ((\max_j s_j)^{-1/2})$ in Theorem \ref{thm1}b.  
So $f_d (\beta_0) = o(1)$, since we allow $s_j \to \infty$.
Then apply our delta theorem, Theorem \ref{thm1}b here
\begin{equation}
r_n [ (w' \hat{\Sigma} w )^{1/2} - (w' \Sigma w)^{1/2}] = o_p (1),\label{ex27}
\end{equation}

Now we see that rate of convergence in risk estimation is faster in (\ref{ex27}), compared to (\ref{ex24})-(\ref{ex25}).

\vspace{2in}

{\bf REFERENCES}

Abadir, K. and J.R. Magnus (2005). Matrix Algebra. Cambridge University Press. Cambridge.

Buhlmann, P. and S. van de Geer (2011). Statistics for High-Dimensional Data. Springer Verlag, Berlin.

Caner, M., E. Ulasan, L. Callot, and O, Onder (2016). "A relaxed approach to estimating large portfolios and gross exposure," arXiv: 1611.07347v1.

Fan, J., Y. Liao, and X. Shi (2015). "Risks of large portfolios," {\it Journal of Econometrics}, 186, 367-387.

Horn, R.A. and C. Johnson (2013). Matrix Analysis. Second Edition. Cambridge University Press, Cambridge.

Van de Geer, S., P. Buhlmann, Y. Ritov, and R. Dezeure (2014). "On asymptotically optimal confidence regions and test for high-dimensional models," {\it The Annals of Statistics,} 42, 1166-1202. 

Van der  Vaart, A.W. (2000). Asymptotic Statistics.  Cambridge University Press, Cambridge.

\newpage

{\bf  \large Appendix}

{\bf Proof of Theorem \ref{thm1}a}. The first part of the theorem shows why the classical proof of delta theorem will not work in high dimensions.  
However, this is not a negative result since it will guide us through the second part which shows us the solution.

{\bf Part 1}. First by differentiability, and define $l(.)$ as a vector function, via p.352 of Abadir and Magnus (2005) $l(.): D \subset R^p \to R^m$
\begin{equation}
\| f (\hat{\beta}) - f (\beta_0) - f_d (\beta_0) [\hat{\beta}- \beta_0] \|_2 = \| l (\hat{\beta} - \beta_0)  \|_2. \label{pt1}
\end{equation}
and 
\begin{equation}
 \frac{\| l(\hat{\beta} - \beta_0) \|_2}{\| \hat{\beta} - \beta_0 \|_2} = o_p (1),\label{pt2}
\end{equation}
where we use Lemma 2.12 of van der Vaart (2000).  Then with $ \| \hat{\beta} - \beta_0 \|_2 = o_p (1)$, (\ref{pt2}) implies
\begin{equation}
\| l (\hat{\beta} - \beta_0) \|_2 = o_p (1).\label{pt3}
\end{equation}
Since we are given $r_n \| \hat{\beta} - \beta_0 \|_2 = O_p (1)$, by (\ref{pt2})(\ref{pt3})
\begin{equation}
r_n \| l (\hat{\beta} - \beta_0)\|_2 = o_p (1).\label{pt4}
\end{equation}
By (\ref{pt1})-(\ref{pt4})
\begin{equation}
\| r_n [ f(\hat{\beta}) - f (\beta_0)] - r_n [ f_d (\beta_0)] [ \hat{\beta} - \beta_0] \|_2 = o_p (1).\label{pt5}
\end{equation}
But this is the same result as in regular delta method. (\ref{pt5}) is mainly a simple extension of Theorem 3.1 in van der Vaart (2000) to Euclidean spaces so far. However the main caveat comes from derivative matrix $f_d (\beta_0)$ which is of dimension $m \times p$. The rate of the matrix plays a role when $p \to \infty$ as $n \to \infty$.  For example, both $r_n [ f_d (\beta_0)] [ \hat{\beta} - \beta_0] $ and $r_n [ f(\hat{\beta}) - f (\beta_0)]$ may be diverging, but $r_n \| \hat{\beta} - \beta_0\| = O_p (1)$. Hence the delta method is not that useful if our interest centers on getting rates for estimators as well as functions of estimators that converge. In the fixed $p$ case, this is not an issue, since the matrix derivative will not affect the rate of convergence at all, as long as this is bounded away from zero, and bounded from above. Note that boundedness assumptions may not be intact when we have $p \to \infty$, as $n \to \infty$. Next part shows how to correct this problem.

{\bf Part 2}.
From differentiability, using p.352 of Abadir and Magnus (2005), or proof of Theorem 3.1 in van der Vaart (2000)

\[ f(\hat{\beta}) - f (\beta_0) = f_d (\beta_0) [ \hat{\beta} - \beta_0] + l (\hat{\beta} - \beta_0).\] 

\noindent Putting the above in Euclidean norm, and using triangle inequality
\begin{eqnarray*}
\| f (\hat{\beta}) - f (\beta_0) \|_2 & = & \| f_d (\beta_0) [ \hat{\beta} - \beta_0] + l (\hat{\beta} - \beta_0) \|_2 \\
& \le & \| f_d (\beta_0) [ \hat{\beta} - \beta_0] \|_2 + \| l (\hat{\beta} - \beta_0) \|_2
\end{eqnarray*}
Next, multiply each side by $r_n$, and use (\ref{pt4})
\begin{equation}
r_n \| f (\hat{\beta}) - f (\beta_0) \|_2 \le r_n \| f_d (\beta_0) [ \hat{\beta} - \beta_0] \|_2  + o_p (1).\label{pt6}
\end{equation}
Then apply  matrix norm inequality in (\ref{1}) to the first term on the right side of (\ref{pt6})
\begin{equation}
r_n \| f_d (\beta_0) [ \hat{\beta} - \beta_0] \|_2 \le r_n \left[ \|| f_d (\beta_0)\||_2 \right]  \left[\| \hat{\beta} - \beta_0 \|_2 \right].\label{pt7}
\end{equation}
Substitute (\ref{pt7}) into (\ref{pt6}) to have 
\begin{equation}
r_n \| f (\hat{\beta}) - f (\beta_0) \|_2 \le r_n \left[ \|| f_d (\beta_0) \||_2 \right] \left[   \| \hat{\beta} - \beta_0 \ \|_2 \right]  + o_p (1).\label{pt8}
\end{equation}
Now divide each side by $\|| f_d (\beta_0) \||_2$,  since $\|| f_d (\beta_0) \||_2 >0$, and $\||f_d (\beta_0)\||_2  \neq o(1)$,
\begin{equation}
\frac{r_n}{ \|| f_d (\beta_0) \||_2} \| f (\hat{\beta}) - f (\beta_0) \|_2 \le r_n \left[   \| \hat{\beta} - \beta_0 \ \|_2 \right]  + o_p (1).\label{pt9}
\end{equation}
By Assumption since $r_n \left[   \| \hat{\beta} - \beta_0 \ \|_2 \right]  = O_p (1)$, we have 
\begin{equation}
\frac{r_n}{ \|| f_d (\beta_0) \||_2} \| f (\hat{\beta}) - f (\beta_0) \|_2 = O_p (1) + o_p (1).\label{pt10}
\end{equation}
So the result is derived, by noting that the new rate of convergence for the function of estimators is
\[ r_n^* =  O \left( \frac{r_n}{ \|| f_d (\beta_0) \||_2} \right).\]
{\bf Q.E.D.}

{\bf Proof of Theorem \ref{thm1}b}.  Here we use (\ref{pt8}), and since $\||f_d (\beta_0)\||_2 = o(1)$, and 
\[ r_n \| \hat{\beta} - \beta_0\|_2 = O_p (1),\]
we have 
\[ r_n \| f(\hat{\beta} ) - f (\beta_0)\|_2 \le o_p (1) + o_p (1) = o_p (1).\]
{\bf Q.E.D.}

\end{document}